\documentclass[12pt]{article}
\usepackage{amsfonts,latexsym,graphicx}
\DeclareMathAlphabet{\mathpzc}{OT1}{pzc}{m}{it}
\input xy
\xyoption{all}
\usepackage{graphics}
\usepackage{amsfonts,latexsym}
\usepackage{graphics,amsmath,epsfig,latexsym,amssymb,cite,graphicx,a4}

\DeclareMathAlphabet{\mathpzc}{OT1}{pzc}{m}{it}
\begin{document}
\title{Premonoidal categories associated with representations of finite groups and their quantum doubles}
\author{
L.D. Wagner, J. Links,
P.S. Isaac,
W.P. Joyce$^*$, 
K. Dancer \\
~~\\
Centre for Mathematical Physics, \\
The University of Queensland, \\  Brisbane, 4072, Australia 
\\
~~\\
$^*$Department of Physics and Astronomy,
\\ The University of Canterbury, \\ 
Private Bag 4800, Christchurch, New Zealand   
}                     
\maketitle
\begin{abstract}
We study the construction of premonoidal categories, where the pentagon relation fails, 
through representations of finite group algebras and their quantum doubles. Both finite group algebras 
and their quantum doubles have a finite number of irreducible representations. We show that in each case there are
at least $2^{n-1}$ inequivalent premonoidal categories of representations, where $n$ is the number of irreducible representations. By construction, for the case of finite group algebras the categories are symmetric whereas for the quantum doubles the categories are braided.   
\end{abstract}

\vfil\eject


\def\a{\alpha}
\def\b{\beta}
\def\d{\delta}
\def\e{\epsilon}
\def\ve{\varepsilon}
\def\g{\gamma}
\def\k{\kappa}
\def\l{\lambda}
\def\o{\omega}
\def\t{\theta}
\def\s{\sigma}
\def\D{\Delta}
\def\L{\Lambda}
\def\X{\bar{X}}
\def\Y{\bar{Y}}
\def\Z{\bar{Z}} 
\def\ch{\check}
\def\f{{\cal F}} 
\def\G{{\cal G}}
\def\hG{{\hat{\cal G}}}
\def\R{{\cal R}}
\def\hR{{\hat{\cal R}}}
\def\tR{{\tilde{\cal R}}}
\def\C{{\mathcal C}}
\def\S{{\mathcal S}}
\def\V{{\mathcal V}}
\def\P{{\bf P}}
\def\T{{\cal T}}
\def\H{{\cal H}}
\def\trho{{\tilde{\rho}}}
\def\tP{{\tilde{\Phi}}}
\def\tT{{\tilde{\cal T}}}
\def\ot{\otimes} 
\def\aa{{\mathpzc{A}}}
\def\qq{{\mathpzc{Q}}}
\def\cmplx{\mathbb{C}}
\def\tr{\mbox{tr}}
\newcommand{\ov}[1]{\ensuremath{\overline{#1}}}
\def\ta{\tau}
\def\s{\sigma}
\def\cchi{\overline{\chi}}

\def\tp{\otimes}

\def\I{{\rm id}}

\def\ev{{\rm ev}}

\def\coev{{\rm coev}}

\def\rar{\rightarrow}

\def\idel{{\mathbf 1}}


\def\beq{\begin{equation}}
\def\eeq{\end{equation}}
\def\bea{\begin{eqnarray}}
\def\eea{\end{eqnarray}}
\def\ba{\begin{array}}
\def\ea{\end{array}}
\def\no{\nonumber}
\def\lt{\left}
\def\rt{\right}
\newcommand{\bq}{\begin{quote}}
\newcommand{\eq}{\end{quote}}

\newtheorem{Theorem}{Theorem}
\newtheorem{Definition}{Definition}
\newtheorem{Proposition}{Proposition}
\newtheorem{Lemma}{Lemma}
\newtheorem{Corollary}{Corollary}

\section{Introduction} 

The goal of this paper is to provide explicit examples of categories where monoidality, 
represented by the pentagon condition \cite{mc,js}, fails to hold, but with a weakened 
notion of coherence remaining. We refer to such categories as premonoidal. Early work on this problem 
was undertaken in \cite{yanofsky} in the context of $n$-categories. The question was pursued in 
\cite{joyce1}, motivated 
by the problem of assigning Bose/Fermi statistics to representations of the Lie algebra $su(n)$. 
This approach was undertaken in a way that generalised the $su(2)$ case. There, the familar notion of assigning odd-dimensional irreducible representations 
as fermionic and odd-dimensional irreducible representations as bosonic leads to a non-trivial 
{\it symmetric} monoidal category of 
finite-dimensional representations. (By symmetric category we mean that the category is braided such that all braiding morphisms square to the identity.) This may be viewed a category-theoretic statement of the spin-statistics theorem.  It was found in \cite{joyce1} that the imposition 
of  Bose/Fermi statistics on the irreducible finite-dimensional representations of 
$su(n)$, $n>2$, leads to a breaking of monoidality. 
Specifically, the pentagon condition was replaced with the following commutative diagram   
\begin{eqnarray}
\xymatrix{
(U\tp V)\tp(W\tp Z)  \ar[dd]_{a_{(U\tp V),W,Z}}                           &
                                                                        & 
(U\tp V)\tp(W\tp Z) \ar[ll]_{q_{U,V,W,Z}}                               \\ 
&&\\
((U\tp V)\tp W)\tp Z                                                    &
                                                                        &
U\tp (V\tp (W\tp Z)) \ar[uu]_{a_{U,V,(W\tp Z)}} \ar[dd]^{\I\tp a_{V,W,Z}} \\
&&\\
(U\tp(V\tp W))\tp Z \ar[uu]^{a_{U,V,W}\tp \I}                           &
                                                                        &
U\tp ((V\tp W)\tp Z) \ar[ll]^{a_{U,(V\tp W),Z}}    
}
\label{qpent}
\end{eqnarray}
where for all objects $U,V,W,Z$ in the category $\C$ the $a_{U,V,W}$ are the associator morphisms and the morphism $q_{U,V,W,Z}$ is not necessarily the identity.
In \cite{joyce2} the coherent diagrams for such 
premonoidal categories were represented by rooted planar binary trees
with levels, and formal primitive operations on these trees. The more detailed concepts underlying braided premonoidal coherence were developed in \cite{joyce3}. 

It is well known that monoidal categories can be systematically constructed  through finite-dimensional representations of quasi-bialgebras, where an additional braiding structure exists if the bialgebra is quasi-triangular 
\cite{drinfeld1,drinfeld2,m,cp,majid}. Symmetric categories arise in the braided case when the bialgebra is cocommutative and the universal $R$-matrix is trivial. Motivated by this approach, an algebraic prescription called {\it twining} was given in \cite{ijl} which deformed the structure of the bialgebras in such a manner that the pentagon condition is broken as represented by (\ref{qpent}). 
The added ingredient used in this construction is a central element $K$ 
which takes an integer eigenvalue $\cchi_\lambda(K)$ on any irreducible finite-dimensional module $V(\lambda)$ with representation $\pi_\lambda$. The central element $K$ determines a generalised notion of Bose/Fermi statistics for all 
irreducible modules $V(\lambda)$ according to whether $\cchi_\lambda(K)$ is even or odd. We may define a map $\S: \C \rightarrow {\mathbb Z}_2$ such that $\S(\lambda)=\cchi_\lambda(K) \mod 2$. In general, letting $\V$ denote the set of irreducible $A$-modules we refer to a map $\S: \V \rightarrow {\mathbb Z}_2$  as a {\it signature} of $A$ provided it maps the trivial module (as defined in Sect. 2) to zero. In other words, the signature is an assignment of Bose-Fermi statistics to the set of $A$-modules such that the trivial module is bosonic. If a signature assigns a bosonic statsitic to all irreducible modules we say that it is trivial. Our convention throughout is to label the trivial module by 0, and when there are $n$ irreducible $A$-modules we omit $\S(0)=0$ and write $\S=(\S(1),\S(2),\dots\S(n-1))$.

Via the above approach it was shown in \cite{ijl} that both the spin-statistics theorem for $su(2)$ and the premonoidal 
$su(3)$ results of \cite{joyce1} could be recovered. We remark that in the general $su(n)$ case discussed in  \cite{joyce1} the construction applies for a limited class of signatures, and the associator morphisms always take the form of a global phase. An open question is given a signature $\S$, can one construct a $K$ which produces $\S$ and in turn provides the means to construct the braided premonoidal category of finite-dimensional representation through twining? 

We will apply the procedure of twining to the case of group algebras of finite groups and their quantum doubles
\cite{dpr,gould}, both of which have the structures of quasi-triangular bialgebras. (Both are in fact examples of Hopf algbras but that additional information will not be needed for our discussions.) 
These cases are amenable to detailed analysis because in both instances the category of finite-dimensional representations has a finite number of irreducible objects. The main result we will establish is that for a given signature, it is always possible to construct a central element which leads to that signature. Thus for the case of finite groups this gives a construction for symmetric premonoidal categories with any signature. Unlike the results discussed earlier for $su(n)$,
the associator morphisms for these categories are not necessarily a global phase operator. The approach extends in a straightforward manner to the case of quantum doubles of the algebras
of finite groups \cite{drinfeld,dpr,gould} (hereafter referred to as the {\it finite group double}). In these instances the premonoidal category of finite-dimensional representations is braided, but no longer symmetric. We provide some detailed analysis in terms of the dihedral groups and their doubles, as their representation theory is well understood and results can be made explicit.

\section{The premonoidal construction from representations of quasi-bialgebras}

It is known that monoidal categories may be constructed using the representations of quasi-bialgebras. 
Here we provide a short review of these algebras, and provide the required theory necessary to deform them in a way that permits the construction of premonoidal categories.

\begin{Definition}\label{quasi-bi} 
A complex quasi-bialgebra $(A,\D,\e,\Phi)$ is an associative algebra $A$
over $\mathbb C$  with unit element $I$, 
equipped with algebra homomorphisms $\e: 
A\rightarrow {\mathbb C}$ (counit), $\D: A\rightarrow A\otimes A$ (coproduct),
and an invertible element $\Phi\in A\otimes A\otimes A$ 
(coassociator) satisfying 
\bea
&& (\I\otimes\D)\D(a)=\Phi^{-1}(\D\otimes \I)\D(a)\Phi,~~
      \forall a\in A,\label{quasi-bi1}\\
&&(\D\otimes \I\otimes \I)\Phi \cdot (\I\otimes \I\otimes\D)\Phi
=(\Phi\otimes I)\cdot(\I\otimes\D\otimes \I)\Phi\cdot (I\otimes
                  \Phi),\label{quasi-bi2} \label{pent}\\
&&m(\e\otimes \I)\D=\I=m(\I\otimes\e)\D,\label{quasi-bi3}\\
&&m(m\tp \I)
(\I\otimes\e\otimes \I)\Phi=I. \label{quasi-bi4} \eea  
If there
exists an invertible element $\R\in A\otimes A$ (universal $R$-matrix) 
such that
\bea
&&\R\D(a)=\D^T(a)\R,~~~~\forall a\in A,\label{dr=rd}\\
&&(\D\otimes
\I)\R=\Phi^{-1}_{231}\R_{13}\Phi_{132}\R_{23}\Phi^{-1}_{123},
   \label{d1r}\\
   &&(\I\otimes \D)\R=\Phi_{312}\R_{13}\Phi^{-1}_{213}\R_{12}\Phi_{123}
      \label{1dr}
      \eea
then $(A,\D,\e,\Phi,R)$ is called a quasi-triangular quasi-bialgebra. 
\end{Definition}

These defining relations are sufficient to ensure that the 
category of finite-dimensional $A$-modules forms a braided monoidal category \cite{drinfeld1,drinfeld2,m,cp,majid}. 
We remark that the co-unit $\e$ 
ensures the existence of a one-dimensional representation, which we will call the {\it trivial} representation and denote by $\pi_0$. This representation is distinguished from other one-dimensional representations which might exist, as will be seen
it must will always be assigned a bosonic statistic.

In order to break monoidality of the category of representations we use twining, as introduced in 
\cite{ijl}.
\begin{Definition} (Twining)
Given a quasi-triangular quasi-bialgebra $(A,\D,\e,\Phi,R)$ which possesses
a central element $K$ taking integer eigenvalues on all irreducible $A$-modules, 
the twined algebra is defined to be $(A,\D,\e,\tP,\tR)$, where 
\bea   
\tR&=&\exp i\pi ({K\otimes K})\cdot\R=\R\cdot\exp(i\pi{K\otimes K}) \no \\
\tP&=&\Phi\cdot\exp(i\pi\kappa) \label{twining} 
\eea 
and
\beq \kappa=
{K\otimes(I\otimes K+K \otimes I-\D(K))}. \label{kappa} \eeq        

The following relations hold: 
\bea
&&(\I\otimes \D)\D(a)=\tP^{-1}(\D\otimes \I)\D(a)\tP~~~~\forall a\in
A,\no\\
&&\tR\D(a)=\D^T(a)\tR,~~~~\forall a\in A,\no\\
&&(\D\otimes
\I)\tR= \tP^{-1}_{231}\tR_{13}\tP_{132}\tR_{23}\tP^{-1}_{123},
   \no \\
   &&(\I\otimes \D)\tR=\tP_{312}\tR_{13}\tP^{-1}_{213}\tR_{12}\tP_{123}[\exp(2i\pi\kappa)]_{123}
   \label{leftR} \eea  \\
We also define 
\beq \xi=(\D\otimes \I\otimes \I)\tP^{-1}\cdot(\tP\otimes I)\cdot(\I\otimes
\D\otimes \I)\tP\cdot(I\otimes\tP)\cdot(\I\otimes \I\otimes \D)\tP^{-1}.
\label{xi} \eeq 
\end{Definition}

Note that in the examples of finite group algebras and their quantum doubles to be discussed below, 
finite-dimensionality of the algebra ensures that the quantities $\exp(i\pi K),\,\exp(i\pi\kappa)$ are well-defined. 
As can be seen from the definition of twining, the twined elements
$\tP$ and $\tR$ depend entirely on the choice of central element $K$.

The application of twining has no effect on equations 
(\ref{quasi-bi1},\ref{quasi-bi3},\ref{dr=rd},\ref{d1r}).
Choosing $K$ such that
\beq
\e(K) = 0,
\label{counitprop}
\eeq
then equation (\ref{quasi-bi4}) holds. We call such central elements {\em admissible}.
Note that equations (\ref{1dr},\ref{leftR}) are not identical.
Following \cite{liu}, we could call $\tR$ a {\em left}
universal quasi $R$-matrix. However, as $K$ is chosen such  
that it takes an integer eigenvalue on every finite-dimensional irreducible 
representation, then all matrix representations of equations (\ref{1dr},\ref{leftR}) 
become equivalent.  The only equation of Definition 2 which differs from its analogue in 
Definition 1 is (\ref{xi}). This is equivalent to (\ref{quasi-bi2}) only when 
$\xi=I\tp I\tp I\tp I$, which is not generally true.


\begin{Definition}\label{bmc}
A premonoidal category is a triple $(\C,\otimes,\aa)$ where $\C$ is a
category, $\otimes:\C\times \C\rightarrow \C$ is a bifunctor and
$\aa:\otimes({\rm id}\times\otimes )\rightarrow\otimes( \otimes\times{\rm id})$ is a
natural isomorphism for associativity. 
A unital pre-monoidal category $\C$ is said to be braided if it is equipped with
a natural commutativity isomorphism $\s_{U,V}:U\tp V \rar V\tp
U$ for all objects $U,V\in \C$ such that the usual triangle and hexagon diagrams commute as for
braided monoidal categories, as does the following diagram

$$
\xymatrix{
(U\tp V)\tp (W\tp Z) \ar[d]_{\s_{(U\tp V),(W\tp Z)}} \ar[rr]^{q_{U,V,W,Z}} &&
(U\tp V)\tp (W\tp Z) \ar[d]^{\s_{(U\tp V),(W\tp Z)}} \\
(W\tp Z)\tp (U\tp V) && (W\tp Z)\tp (U\tp V) \ar[ll]_{q_{W,Z,U,V}}
}
.$$
In general, if $\s_{U,V}\circ\s_{V,U}=\I_{V\tp U}$ for all objects $U,V\in
\C$, $\C$ is said to be {\em symmetric} (also 
referred to as {\em tensor} in the literature). 
\end{Definition}

Consider a quasi-triangular bialgebra $A$ with irreducible $A$-module $V(\lambda)$ and corresponding representation $\pi_\lambda$. Setting 
\bea a_{\lambda,\mu,\nu}&=&(\pi_\lambda\otimes\pi_\mu\otimes\pi_\nu)\tP, \label{a}\\ 
\sigma_{\lambda,\mu}&=&P(\pi_\lambda\otimes\pi_\mu)\tR, \label{sig} \\   
q_{\lambda,\mu,\nu,\rho}&=&(\pi_\lambda\otimes\pi_\mu\otimes\pi_\nu\otimes\pi_\rho)\xi, \label{q} 
 \eea
where $P_{\lambda,\mu}(V(\lambda)\otimes V(\mu))=V(\mu) \otimes V(\lambda)$ is the flip map, 
it was established in \cite{ijl}
that the category of finite-dimensional representations of the twined algebra is a braided,
premonoidal category.      
We emphasise that the imposition $\e(K)=0$ ensures that an
essential property is satisfied (see for example, \cite{cp} page 530, a comment in
the proof of proposition 16.1.2).


\section{Twining the group algebras of finite groups}

\noindent The complex algebra $\cmplx[G]$ of a finite group $G$ provides an
example of a complex quasi-triangular quasi-bialgebra with coproduct and counit 
respectively given by
$$
\D(g) = g\tp g,\ \ \ \e(g)=1,\ \ \ \forall g\in G
$$
and the trivial coassociator and $R$-matrix are
\beq
\Phi = e\tp e\tp e, \ \ \ R = e\tp e
\eeq
where $e$ is the identity of $G$. 

It is well known that we can partition a finite group $G$ into
conjugacy classes. 
The character of a representation $\pi_\lambda$ is the function
$\chi_\lambda:G\rar\cmplx$ defined by
$$
\chi_\lambda(g) = \tr(\pi_\lambda(g)).
$$
\noindent Using the cyclic properties of traces it can be shown that $\chi_{\lambda}$ is a class function. The set of irreducible characters (those arising from irreducible representations) forms a basis for the vector space $\hat{G}$ of class functions. Consequently, the number of irreducible representations of $G$ is precisely the number of conjugacy classes.
Using this it can be shown that the elements 
\begin{equation}
E_\lambda= \frac{d_\lambda}{|G|}\sum_{g\in\,G} \chi_\lambda(g^{-1})g 
\label{cidem}
\end{equation} 
form a basis for the centre of $\cmplx[G]$ \cite{cr}, where $d_\lambda$ denotes the dimension of $\pi_\lambda$ and $|G|$ is the order of $G$. Moreover, these elements are orthogonal idempotents giving a resolution of the identity:
\begin{eqnarray} 
E_\lambda E_\mu &=&\delta_{\lambda,\mu} E_\mu \label{ortho}\\
\sum_{\lambda} E_\lambda &=&e.  
\label{resol}
\end{eqnarray}  
Their eigenvalues on the irreducible modules are simply 
$$\cchi_\lambda(E_\mu)=\delta_{\lambda,\mu} $$
from which it follows that the action of the co-unit is 
$$\e(E_\lambda)=\delta_{0,\lambda}$$ 
where $\pi_0\equiv \e$. 
Given a signature $\S$ we construct the central element 
$$K_\S= \sum_{\lambda} \S(\lambda) E_\lambda $$  
satisfying 
$$K_\S^2=K_\S,~~~~~~~\e(K_\S)=0. $$ 
As $K_\S$ is admissible, we have  
\begin{Proposition}
Given a signature $S$ of $\cmplx[G]$, the category of $\cmplx[G]$-modules is a symmetric premonoidal category with 
$a_{\lambda,\mu,\nu},\, q_{\lambda,\mu,\nu,\rho} ,\,\sigma_{\lambda,\mu}$ given by (\ref{a},\ref{sig},\ref{q}) 
respectively where 
\bea \tP&=&
( e^{\otimes 3}  -2 K_\S\otimes K_S\otimes e)
(e^{\otimes 3} -2 K_\S\otimes e \otimes K_\S)
(e^{\otimes 3} -2 K_\S \otimes \Delta(K_\S)) \no \\
 \xi&=& ( e^{\otimes 4}  -2 e\otimes K_\S\otimes e \otimes K_\S ) 
( e^{\otimes 4}  -2 K_\S\otimes e \otimes K_\S \otimes e) 
( e^{\otimes 4}  -2 e\otimes K_\S \otimes K_\S \otimes e   ) \no \\
&& ~~~~\times (e^{\otimes 4}  -2 K_\S\otimes e \otimes e\otimes K_\S) 
(e^{\otimes 4}  -2 e\otimes K_\S \otimes \Delta(K_\S)  ) 
(e^{\otimes 4}  -2 K_\S \otimes e \otimes \Delta(K_\S)  )  \no \\
&&~~~~\times (e^{\otimes 4}  -2 \Delta(K_\S) \otimes e \otimes K_\S  ) 
(e^{\otimes 4}  -2 \Delta(K_\S) \otimes K_\S \otimes \otimes e  ) 
(e^{\otimes 4}  -2 \Delta(K_\S)\otimes \Delta(K_\S)) 
\no \\
 \tR&=&(e\otimes e-2K_\S\otimes K_\S) \no \eea
and
$$\Delta(K_\S)=\sum_{\lambda} \frac{\S(\lambda)d_\lambda}{|G|}\sum_{g\in\,G} \chi_\lambda(g^{-1})g \otimes g. $$
\end{Proposition}

In the above construction there are $2^{n-1}$ choices for the signature, where $n$ is the number of irreducible representations. Hence there are $2^{n-1}$ inequivalent symmetric premonoidal categories of representations, including the 
case with trivial signature. 

\section{The quantum double algebra $D(G)$}

The above results for group algebras extend to the case of their quantum doubles. The quantum double is a 
construction which embeds any Hopf algebra to be embedded in a quasi-triangular Hopf algebra \cite{drinfeld}. 
First we give a brief survey of the quantum double $D(G)$ of a finite group $G$ 
\cite{dpr,gould}.

Let $\cmplx[G]^*$ denote the dual space of $\cmplx[G]$ , so $\cmplx[G]^* = \{f|\,f:G 
\rightarrow \mathbb{C}\}$.  Explicitly we define 

$$ g^*(h) = \delta(g,h) \quad \forall g,h \in G.$$ 

\noindent Then  $\cmplx[G]^*$ is the algebra of the dual elements $g^*$ with multiplication

$$g^* h^* = \delta(g,h) h^*.$$
The quantum double  $D(G)$ is a $|G|^2$-dimensional algebra spanned by the free products

$$ g h^*, \quad g, h \in G,$$

\noindent where the elements $h^* g$ are calculated using

$$ h^* g = g (g^{-1}hg)^*.$$
Then $D(G)$ is a quasi-triangular Hopf algebra with coproduct $\ov{\Delta}$ and counit $\ov{\varepsilon}$ given by:

\begin{align*}
&\ov{\Delta}(gh^*) = \sum_{k \in G} g (k^{-1}h)^* \otimes gk^* = \sum_{k \in G}
  gk^* \otimes g(hk^{-1})^*, \\
&\ov{\ve}(gh^*) = \delta(h,e).
\end{align*}
Note that we identify $g \ve$ with $g$ and $e g^*$ with $g^*$ for all $g \in G$.
The universal $R$-matrix is given by 

$$ R = \sum_{g \in G} g \otimes g^* $$

\noindent which can easily be shown to satisfy the defining relations (2-8) with $\phi=e\otimes e \otimes e$. 

In analogy with the finite group case, the character of a representation $\pi_\lambda$ is the function
$\chi_\lambda:D(G)\rar\cmplx$ defined by 
$$
\chi_\lambda(g) = \tr(\pi_\lambda(g)).
$$
{}From \cite{gould} we have the finite group double analogues of (\ref{cidem})  
\bea 
E_\lambda= \frac{d_\lambda}{|G|}\sum_{g,h\in\,G} \chi_\lambda(g^{-1}h^*)gh^* 
\eea 
which satisfy the relations (\ref{ortho},\ref{resol}). We remark that 
$\chi_\lambda(g^{-1}h^*)=0$ unless $gh=hg$ \cite{gould}. 
Proceeding as in the finite group algebra case,
given a signature $\S$ we construct the central element 
$$K_\S= \sum_{\lambda} \S(\lambda) E_\lambda $$  
satisfying 
$$K_\S^2=K_\S,~~~~~~~\e(K_\S)=0. $$ 
We then have  

\begin{Proposition}
Given a signature $S$ of $D(G)$, the category of $D(G)$-modules is a braided premonoidal category with 
$a_{\lambda,\mu,\nu},\, q_{\lambda,\mu,\nu,\rho} ,\,\sigma_{\lambda,\mu}$ given by (\ref{a},\ref{sig},\ref{q}) respectively
where 
\bea \tP&=&
( e^{\otimes 3}  -2 K_\S\otimes K_S\otimes e)
(e^{\otimes 3} -2 K_\S\otimes e \otimes K_\S)
(e^{\otimes 3} -2 K_\S \otimes \Delta(K_\S)) \no \\
 \xi&=& ( e^{\otimes 4}  -2 e\otimes K_\S\otimes e \otimes K_\S ) 
( e^{\otimes 4}  -2 K_\S\otimes e \otimes K_\S \otimes e) 
( e^{\otimes 4}  -2 e\otimes K_\S \otimes K_\S \otimes e   ) \no \\
&& ~~~~\times (e^{\otimes 4}  -2 K_\S\otimes e \otimes e\otimes K_\S) 
(e^{\otimes 4}  -2 e\otimes K_\S \otimes \Delta(K_\S)  ) 
(e^{\otimes 4}  -2 K_\S \otimes e \otimes \Delta(K_\S)  )  \no \\
&&~~~~\times (e^{\otimes 4}  -2 \Delta(K_\S) \otimes e \otimes K_\S  ) 
(e^{\otimes 4}  -2 \Delta(K_\S) \otimes K_\S \otimes \otimes e  ) 
(e^{\otimes 4}  -2 \Delta(K_\S)\otimes \Delta(K_\S)) 
\no \\
 \tR&=&(e\otimes e-2K_\S\otimes K_\S)\sum_{g\in\,G} g\otimes g^* \no \eea
and
$$\Delta(K_\S)=\sum_{\lambda} \frac{\S(\lambda)d_\lambda}{|G|}\sum_{g,h\in\,G} 
\chi_\lambda(g^{-1}h^*)\sum_{k\in\,G} g(k^{-1}h)^* \otimes gk^*. $$
\end{Proposition}

In the above construction there are $2^{n-1}$ choices for the signature, where $n$ is the number of irreducible representations. Hence there are $2^{n-1}$ inequivalent braided premonoidal categories of representations, including the 
case with trivial signature. 

\section{Examples} 

Here we illustrate the theory with worked examples. 
The representation theory of the general dihedral groups $D_n$ is 
known which allows us to make the above results explicit.
The dihedral group $D_n$  has two generators $\sigma, \tau$ satisfying:

$$ \sigma^n = e,\; \tau^2 = e,\; \tau \sigma = \sigma^{n-1} \tau .$$
When $n$ is odd, there are $({n+3})/{2}$ conjugacy classes divided into three families, given by:

\begin{align*}
&\{ e \}, \\
&\{ \sigma^k, \sigma^{-k} \} \quad \text{for } 1 \leq k \leq \frac{n-1}{2}, \\
&\{ \sigma^i \tau, \, 0 \leq i \leq n-1 \}.\notag
\end{align*}
There are $({n+3})/{2}$ irreducible representations, 
two of which are one-dimensional and the remaining $({n-1})/{2}$ which are two-dimensional. They are given by:

$$ \pi_\pm (\sigma) =1,\quad \pi_\pm (\tau) = \pm 1  $$ 

\noindent and

\begin{equation*}
\pi_k (\sigma) = \begin{bmatrix}
                    \omega^k & 0  \\ 0 & \omega^{-k}
                 \end{bmatrix}, \quad
\pi_k (\tau) = \begin{bmatrix}
                 0 & 1 \\ 1 & 0 
               \end{bmatrix}, \quad
\omega = \exp\left({\frac{2 \pi i}{n}}\right), \; 1 \leq k \leq \frac{n-1}{2}.
\end{equation*}

When $n$ is even, there are $({n}+6)/{2}$ conjugacy classes divided into five families, given by:

\begin{align*}
&\{e\}, \notag \\
&\{\sigma^{{n}/{2}} \}, \notag \\
&\{ \sigma^k, \sigma^{-k} \} \quad \text{for } 1 \leq k \leq {(n-2)}/{2}, \\
&\{ \sigma^{2j} \tau, \; 0 \leq j \leq {(n-2)}/{2} \}, \notag \\
&\{ \sigma^{(2j+1)} \tau, \; 0 \leq j \leq {(n-2)}/{2} \}. \notag 
\end{align*}

\noindent The $({n+6})/{2}$ irreps consist of 4 one-dimensional irreps and ${(n-2)}/{2}$ two-dimensional 
irreducible reprentations.  They are given by:

$$ \pi(\sigma) = (-1)^a, \quad \pi (\tau) = (-1)^b \quad \text{for } \; a,b \in \{ 0,1 \} $$

\noindent and 

\begin{equation*}
\pi_k (\sigma) = \begin{bmatrix}
                    \omega^k & 0  \\ 0 & \omega^{-k}
                 \end{bmatrix}, \quad
\pi_k (\tau) = \begin{bmatrix}
                 0 & 1 \\ 1 & 0 
               \end{bmatrix}, \quad
\omega = \exp\left({\frac{2\pi i}{n}}\right),\; 1 \leq k \leq \frac{n-2}{2}.
\end{equation*}

The above data is sufficient to explicitly determine the central idempotents $E_\lambda$, and in turn the co-associator and braiding isomorphisms. Below we give the results for the case of $D_3$.

\subsection{Twining the dihedral group $D_{3}$}

\noindent The group $D_{3}$ is of order $6$ consisting of the elements
\begin{equation}
D_{3}=\{ e,\s, \s^{2},\ta, \s \ta, \s^{2}\ta \}.
\end{equation}
\noindent The conjugacy classes for $D_{3}$ are $C_{0}=\{ e \}, \; C_{1}=\{\s ,\s^2 \} ~{\rm and }~ 
C_{2}=\{\ta ,\s \ta , \s^2  \ta \}$.

\begin{table}[h!] \label{chard3}
\begin{center}
\begin{tabular}{|c|ccc|}
\hline
$D_{3}$ & $C_{0}$ & $C_{1}$ & $C_{2}$\\ \hline
$\pi_{0}$ & $1$ & $1$ & $1$\\
$\pi_{1}$ & $1$ & $1$ & $-1$\\
$\pi_{2}$ & $2$ & $-1$ & $0$\\
 \hline
\end{tabular}
\end{center}
\caption{Character table for $D_{3}$} 
\end{table}
 
\noindent Using the character table, we can explicitly 
construct the central operators $E_j$:
\bea 
E_0&=&\frac{1}{6}\left[ e+\sigma+\sigma^2 +\tau+\sigma\tau+\sigma^2\tau \right] \\
E_1&=&\frac{1}{6}\left[ e+\sigma+\sigma^2 -(\tau+\sigma\tau+\sigma^2\tau) \right] \\
E_2&=&\frac{1}{3}\left[2e-(\sigma+\sigma^2)\right] \eea
and for signature $\S=(\S(1),\,\S(2))$
\bea 
K_\S&=&\frac{1}{6}\left( [4\S(1)+\S(2)]e+[\S(2)-2\S(1)][\sigma+\sigma^2]\right.  \no \\
&&~~~~~~\left.-\S(2)[\tau+\sigma\tau+\sigma^2\tau]  \right). 
\eea 
{}From Proposition 1 we can explicitly determine the co-associator isomorphisms and braiding isomorphisms. The latter are simply given by 
\begin{eqnarray*}
\sigma_{\lambda,\mu}=(-1)^{\S(\lambda)\S(\mu)} P_{\lambda,\mu} 
\end{eqnarray*}
for any signature. 
The co-associator isomorphisms are listed below.
~~\\
$\S=(0,1)$:
$$\begin{array}{rrrr}
a_{0,0,0}= 1\otimes 1 \otimes 1 & a_{0,0,1}=  1\otimes 1 \otimes 1& 
a_{0,0,2}=  1\otimes 1 \otimes I_2&  a_{0,1,0}= 1\otimes 1 \otimes 1 \\
a_{0,1,1}= 1\otimes 1 \otimes 1& a_{0,1,2}= 1\otimes 1 \otimes I_2 & 
a_{0,2,0}= 1\otimes I_2 \otimes 1& a_{0,2,1}=  1\otimes I_2 \otimes 1\\
a_{0,2,2}=  1\otimes I_2 \otimes I_2& a_{1,0,0}= 1\otimes 1 \otimes 1 & 
a_{1,0,1}= 1\otimes 1 \otimes 1& a_{1,0,2}= 1\otimes 1 \otimes I_2\\ 
a_{1,1,0}= 1\otimes 1 \otimes 1& a_{1,1,1}= 1\otimes 1 \otimes 1 & 
a_{1,1,2}=  1\otimes 1 \otimes I_2& a_{1,2,0}= 1\otimes I_2 \otimes 1 \\
a_{1,2,1}= 1\otimes I_2 \otimes 1& a_{1,2,2}= 1\otimes I_2 \otimes I_2& 
a_{2,0,0}= I_2\otimes 1 \otimes 1& a_{2,0,1}= I_2\otimes 1 \otimes 1 \\
a_{2,0,2}= I_2\otimes 1 \otimes I_2 & a_{2,1,0}=  I_2\otimes 1 \otimes 1& 
a_{2,1,1}= I_2\otimes 1 \otimes 1& a_{2,1,2}=  I_2\otimes 1 \otimes I_2\\ 
a_{2,2,0}= I_2\otimes I_2 \otimes 1& a_{2,2,1}=  I_2\otimes I_2 \otimes 1& 
a_{2,2,2}=  -I_2\otimes Q&    
\end{array}$$

$\S=(1,0)$:
$$\begin{array}{rrrrrr}
a_{0,0,0}= 1\otimes 1 \otimes 1 & a_{0,0,1}=  1\otimes 1 \otimes 1& 
a_{0,0,2}=  1\otimes 1 \otimes I_2&  a_{0,1,0}= 1\otimes 1 \otimes 1 \\
a_{0,1,1}= 1\otimes 1 \otimes 1& a_{0,1,2}= 1\otimes 1 \otimes I_2 & 
a_{0,2,0}= 1\otimes I_2 \otimes 1& a_{0,2,1}=  1\otimes I_2 \otimes 1\\
a_{0,2,2}=  1\otimes I_2 \otimes I_2& a_{1,0,0}= 1\otimes 1 \otimes 1 & 
a_{1,0,1}= 1\otimes 1 \otimes 1& a_{1,0,2}= 1\otimes 1 \otimes I_2\\ 
a_{1,1,0}= 1\otimes 1 \otimes 1& a_{1,1,1}= 1\otimes 1 \otimes 1 & 
a_{1,1,2}=  -1\otimes 1 \otimes I_2& a_{1,2,0}= 1\otimes I_2 \otimes 1 \\
a_{1,2,1}= -1\otimes I_2 \otimes 1& a_{1,2,2}= 1\otimes N& 
a_{2,0,0}= I_2\otimes 1 \otimes 1& a_{2,0,1}= I_2\otimes 1 \otimes 1 \\
a_{2,0,2}= I_2\otimes 1 \otimes I_2 & a_{2,1,0}=  I_2\otimes 1 \otimes 1& 
a_{2,1,1}= I_2\otimes 1 \otimes 1& a_{2,1,2}=  I_2\otimes 1 \otimes I_2\\ 
a_{2,2,0}= I_2\otimes I_2 \otimes 1& a_{2,2,1}=  I_2\otimes I_2 \otimes 1& 
a_{2,2,2}=  I_2\otimes I_2 \otimes I_2&    
\end{array}$$

$\S=(1,1)$:
$$\begin{array}{rrrrrr}
a_{0,0,0}= 1\otimes 1 \otimes 1 & a_{0,0,1}=  1\otimes 1 \otimes 1& 
a_{0,0,2}=  1\otimes 1 \otimes I_2&  a_{0,1,0}= 1\otimes 1 \otimes 1 \\
a_{0,1,1}= 1\otimes 1 \otimes 1& a_{0,1,2}= 1\otimes 1 \otimes I_2 & 
a_{0,2,0}= 1\otimes I_2 \otimes 1& a_{0,2,1}=  1\otimes I_2 \otimes 1\\
a_{0,2,2}=  1\otimes I_2 \otimes I_2& a_{1,0,0}= 1\otimes 1 \otimes 1 & 
a_{1,0,1}= 1\otimes 1 \otimes 1& a_{1,0,2}= 1\otimes 1 \otimes I_2\\ 
a_{1,1,0}= 1\otimes 1 \otimes 1& a_{1,1,1}= 1\otimes 1 \otimes 1 & 
a_{1,1,2}=  -1\otimes 1 \otimes I_2& a_{1,2,0}= 1\otimes I_2 \otimes 1 \\
a_{1,2,1}= -1\otimes I_2 \otimes 1& a_{1,2,2}= -1\otimes P& 
a_{2,0,0}= I_2\otimes 1 \otimes 1& a_{2,0,1}= I_2\otimes 1 \otimes 1 \\
a_{2,0,2}= I_2\otimes 1 \otimes I_2 & a_{2,1,0}=  I_2\otimes 1 \otimes 1& 
a_{2,1,1}= I_2\otimes 1 \otimes 1& a_{2,1,2}=  -I_2\otimes 1 \otimes I_2\\ 
a_{2,2,0}= I_2\otimes I_2 \otimes 1& a_{2,2,1}= - I_2\otimes I_2 \otimes 1& 
a_{2,2,2}=  -I_2\otimes P&    
\end{array}$$
where $I_2$ is the $2\times 2$ identity matrix. Here  
\bea 
N&=& e^1_1\otimes e^1_1 -  e^1_2\otimes e^2_1 -   e^2_1\otimes e^2_1 +   e^2_2\otimes e^2_2  \\
P&=& e^1_1\otimes e^1_1 +  e^1_2\otimes e^2_1 +   e^2_1\otimes e^2_1 +   e^2_2\otimes e^2_2  \\
Q&=& e^1_1\otimes e^1_1 -  e^1_1\otimes e^2_2 -   e^2_2\otimes e^1_1 +   e^2_2\otimes e^2_2  
\eea 
where $e^j_k$ is the matrix with 1 in the $(j,k)$ position and zeroes elsewhere. Similarly the $q_{\lambda,\mu,\nu,\rho}$ can be explicitly determined. 

\subsection{Twining the quantum double $D(D_3)$ of the dihedral group $D_3$.}

Just as the representation theory of the dihedral groups is completely understood, the representation theory of their 
quantum doubles is also known \cite{dil}. However the results cannot be expressed as compactly as the $D_n$ case as given in 
Sect. 5.  
As $D_3$ is the simplest non-abelian finite group, we instead restrict to $D(D_3)$ which provides the simplest example of a finite group double. Below we give a full description of the irreducible representations in terms of the generators.  

~~\\
\noindent \underline{One-dimensional irreducible representations}

$$\sigma=1, \quad \tau = \pm 1, \quad g^* = \delta(g, e)$$

\noindent \underline{Two-dimensional irreducible representations}

\begin{equation*}
 \sigma = \begin{bmatrix}
                    \exp(2\pi i/3) & 0  \\ 0 & \exp(4\pi i/3)
                 \end{bmatrix}, \quad
\tau = \begin{bmatrix}
                 0 & 1 \\ 1 & 0 
               \end{bmatrix}, \quad 
g^* = \delta(g, e) I_2,
\end{equation*} 

\begin{equation*}
 \sigma = \begin{bmatrix}
                    \exp(2k\pi i/3) & 0  \\ 0 & \exp(4k\pi i/3)
                 \end{bmatrix}, \quad
\tau = \begin{bmatrix}
                 0 & 1 \\ 1 & 0 
               \end{bmatrix}, \quad 
(\sigma)^* = \begin{bmatrix}
                 1 & 0 \\ 0 & 0
               \end{bmatrix}, \quad
(\sigma^{-1})^* = \begin{bmatrix}
                 0 & 0 \\ 0 & 1
               \end{bmatrix}
\end{equation*} 

\noindent and $g^*=0$ otherwise where $0 \leq k < 3$.

\

\noindent \underline{Three-dimensional irreducible representations}

~~\\

\begin{equation*}
\sigma = \begin{bmatrix}
           0 & 1 & 0 \\ 0 & 0 & 1 \\ 1 & 0 & 0 
	 \end{bmatrix}, \quad
\tau = \pm \begin{bmatrix}
             1 & 0 & 0 \\ 0 & 0 & 1 \\ 0 & 1 & 0
	   \end{bmatrix}
\end{equation*}
and 
$$
(\sigma^i)^* = 0,\quad (\sigma^i \tau)^* = E^{i+1}_{i+1}, \quad 0 \leq i < 3.
$$

{}From the above results we can construct the eight central idempotents which span the centre of $D(D_3)$:
\begin{eqnarray*}
E_0 &=& \frac{1}{6}\left[ e+\sigma+\sigma^2 +\tau+\sigma\tau+\sigma^2\tau \right] \\
E_1&=&\frac{1}{6}\left[ e+\sigma+\sigma^2 -(\tau+\sigma\tau+\sigma^2\tau) \right] \\
E_2&=&\frac{1}{3}\left[2e-(\sigma+\sigma^2)\right] \\
E_3 &=& \frac{1}{3}\left[\sigma \sigma ^{*} + \sigma^{2} \sigma^{*} + \sigma^{*} + \sigma (\sigma^{-1})^{*} 
+\sigma^{2} (\sigma^{-1})^{*} +  (\sigma^{-1})^{*} \right]\\
E_4 &=& \frac{1}{3}\left[\sigma^{*} + \exp(2\pi i/3) \sigma \sigma ^{*} + \exp(4\pi i/3)\sigma^{2} \sigma^{*} + \exp(4\pi i/3)\sigma (\sigma^{2})^{*} \right.\\
 &&~~~~ \left. +\exp(2\pi i /3)\sigma^{2} (\sigma^{2})^{*} +(\sigma^{2})^{*} \right]\\
E_5 &=& \frac{1}{3}\left[\sigma^{*} + \exp(4\pi i/3) \sigma \sigma ^{*} + \exp(2\pi i/3)\sigma^{2} \sigma^{*} + \exp(2\pi i /3)\sigma (\sigma^{2})^{*} \right.\\
 &&~~~~ \left. +\exp(4\pi i/3)\sigma^{2} (\sigma^{2})^{*} +(\sigma^{2})^{*} \right]\\
E_6 &=& \frac{1}{2}\left[(\tau)^{*} + \tau(\tau)^{*}+(\sigma\tau)^{*}  +\sigma\tau (\sigma\tau)^{*}
+(\sigma^{2}\tau)^{*}  + \sigma^{2}\tau(\sigma^{2}\tau)^{*}\right] \\
E_7 &=&  \frac{1}{2}\left[(\tau)^{*} - \tau(\tau)^{*}+(\sigma\tau)^{*}  -\sigma\tau(\sigma\tau)^{*} 
+(\sigma^{2}\tau)^{*} -\sigma^{2}\tau (\sigma^{2}\tau)^{*}\right]
\end{eqnarray*}

 There are $2^7-1=127 $ possible non-trivial signatures, for each of which there are $8^3=512 $  co-associator isomorphisms $a_{\lambda,\mu,\nu}$ and $8^2=64$ braiding isomorphisms $\sigma_{\lambda,\mu}$. Consequently we do not give the explicit results here. However for a given 
signature and set of representations labels $\lambda,\mu,\nu$ it is straightforward, if tedious, to compute 
$a_{\lambda,\mu,\nu}$ 
and $\sigma_{\lambda,\mu}$ from the above data using Proposition 2, similar to the previous example. 

\begin{flushleft}
{\bf Acknowledgements} 
\end{flushleft}
We thank Mark Gould and Tel Lekatsas for their helpful comments. LDW gratefully acknowledges the support from a Graduate School Research Travel Award (GSRTA) granted by the University of Queensland. Furthermore, LDW thanks the Department of Physics and Astronomy, University of Canterbury for its kind hospitality during October and November 2005. WPJ would like to acknowledge the support and hospitality of the University of Queensland during a visit in July 2005. This work has also been supported by the Australian Research Council.

\end{document}